\documentclass{article}

\usepackage[
  lang = american, 
 paper = hardcover, openaccess,
]{ems-ecm}

\numberwithin{equation}{section}

\usepackage{cancel}

\newcommand{\bR}{\mathbb{R}}
\newcommand{\bQ}{\mathbb{Q}}
\newcommand{\bZ}{\mathbb{Z}}
\newcommand{\bC}{\mathbb{C}}
\newcommand{\bF}{\mathbb{F}}
\newcommand{\cB}{\mathcal{B}}
\newcommand{\cH}{\mathcal{H}}
\providecommand{\abs}[1]{\left\lvert#1\right\rvert}
\DeclareMathOperator{\GL}{GL}	
\DeclareMathOperator{\SL}{SL}
\DeclareMathOperator{\SO}{SO}
\DeclareMathOperator{\Sp}{Sp}
\DeclareMathOperator{\Id}{Id}				
\DeclareMathOperator{\Aut}{Aut}
\DeclareMathOperator{\Ind}{Ind}

\DeclareMathOperator{\Rep}{Rep}
\DeclareMathOperator{\triv}{triv}
\DeclareMathOperator{\aff}{aff}

\theoremstyle{plain}
\newtheorem{theorem}{Theorem}[section]
\newtheorem{vaguethm}{Vague Theorem}[section] 

\newtheorem{fact}[theorem]{Fact}
\newtheorem*{goal}{Long term goal}
\newtheorem{problem}{Problem}
\theoremstyle{definition}
\newtheorem{definition}[theorem]{Definition}

\begin{document}

\title{An introduction to representations of $p$-adic groups}

\emsauthor{1}{
	\givenname{Jessica}
	\surname{Fintzen}
	\mrid{}
	\orcid{}}{Jessica Fintzen}

\Emsaffil{1}{
	\department{Mathematisches Institut}
	\organisation{University of Bonn}
	\rorid{}
	\address{Endenicher Allee 60}
	\zip{53115}
	\city{Bonn}
	\country{Germany}
	\affemail{fintzen@math.uni-bonn.de}}

\classification[11F27, 22E35, 20C08, 20C20]{22E50}


\begin{abstract}
An explicit understanding of the (category of all smooth, complex) representations of $p$-adic groups provides an important tool not just within representation theory. It also has applications to number theory and other areas, and, in particular, it enables progress on various different forms of the Langlands program. In this write-up of the author's ECM 2024 colloquium-style talk, we will introduce $p$-adic groups and explain how the category of representations of $p$-adic groups decomposes into subcategories, called Bernstein blocks. We also provide an overview of what we know about the structure of these Bernstein blocks including a sketch of recent results of the author with Jeffrey Adler, Manish Mishra and Kazuma Ohara that allow to reduce a lot of problems about the (category of) representations of $p$-adic groups to problems about representations of finite groups of Lie type, where answers are often already known or are at least easier to achieve. Moreover, we provide an overview of what is known about the construction of supercuspidal representations, which are the building blocks of all smooth representations and whose construction is also the key to obtain the above results about the structure of the whole category of smooth representations. We will, in particular, focus on recent advances which include the work of the author mentioned in the EMS prize citation as well as a hint towards her recent joint work with David Schwein. 
\end{abstract}

\maketitle

This article is a written version of the talk that the author gave during the 9th ECM.
For more extended and precise survey articles see \cite{Fi-CDM} (for a general math audience), \cite{Fi-IHES} (for a more specialized audience, including graduate students working in number theory and representation theory) and \cite{Fi-ICM}.

\section{$p$-adic numbers and $p$-adic groups}
In order to talk about the representation theory of $p$-adic groups, we first explain what $p$-adic groups are, for which we should first explain what ``$p$-adic'' even means. Throughout this article, $p$ will denote a prime number. Once we have chosen our favorite prime number $p$, we can define the $p$-adic absolute value of the integers $\bZ$ that measure how often an integer is divisible by $p$.

\begin{definition}
	Let $r$ be a non-zero integer coprime to $p$, and let $s$ be a non-negative integer. The \textit{$p$-adic absolute value $\abs{p^s \cdot r}_p $} of the integer $p^s \cdot r$ is defined to be 
	\[ \abs{p^s \cdot r}_p = \left(\frac{1}{p}\right)^{s} . \]
\end{definition}

We observe that the more often an integer is divisible by $p$, the smaller is its $p$-adic absolute value. 
\begin{definition}
The \textit{$p$-adic integers $\bZ_p$} are defined to be the completion of the integers $\bZ$ by the $p$-adic absolute value. 
\end{definition}
This means that a $p$-adic integer is of the form
\begin{equation*}
	a_0+a_1\cdot p +a_2 \cdot p^2 + a_3 \cdot p^3 + \hdots \text{ for some integers } a_i \, (0 \leq a_i < p),
\end{equation*}
because $\abs{p^n}_p= \left(\frac{1}{p}\right)^{n}$ goes to zero as $n$ goes to infinity.

\begin{definition}
		The \textit{$p$-adic numbers $\bQ_p$} are the fraction field of the $p$-adic integers $\bZ_p$, or, equivalently, the completion of the rational numbers $\bQ$ with respect to $p$-adic absolute value (that can be extended to the rational numbers).
\end{definition}
Similarly to the $p$-adic integers, a $p$-adic number can be represented as
\[ a_{-n} \cdot p^{-n} + \hdots + a_0  +a_1\cdot p +a_2 \cdot p^2 + \hdots \text{ with }  a_i \in \{0, \hdots, p-1\} . \]

\begin{figure}[t]
	\includegraphics[width=\textwidth]{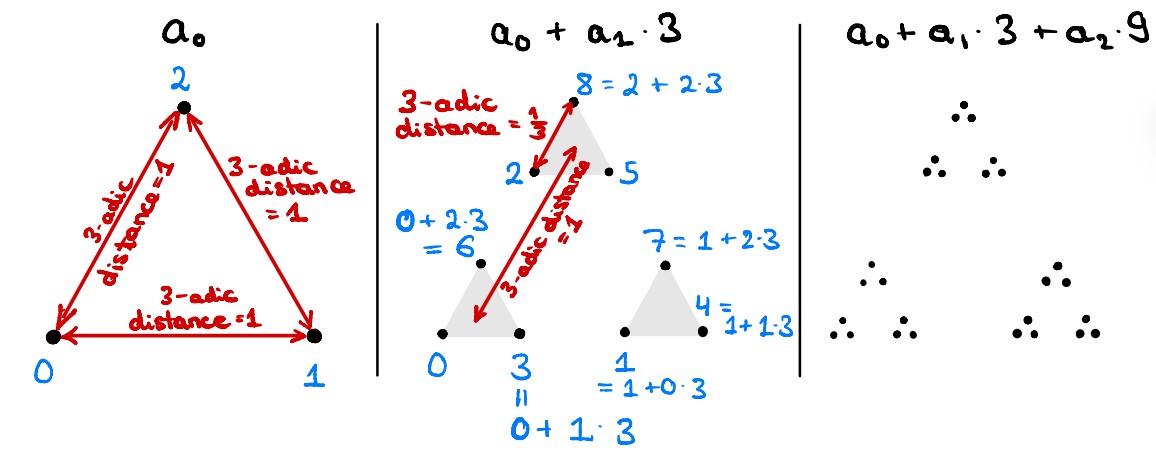}
	\caption{A visualization of the 3-adic integers based on \cite{Cuoco}.}
	\label{Figure-3-adic}
\end{figure}

While the $p$-adic numbers are, like the real numbers $\bR$, a completion of the rational numbers, the topological space of the $p$-adic numbers arising from the $p$-adic absolute value is very different from the topology on the real numbers. The $p$-adic numbers turn out to be totally disconnected. To illustrate the topology on the 3-adic integers as an example, note that the integers 0, 1 and 2 all have 3-adic distance 1 from each other. The integers 0, 3, and 6 all have 3-adic distance $\frac{1}{3}$ from each other, and similarly the integers 1, $1+3=4$, $1+2\cdot 3=7$ all have 3-adic distance $\frac{1}{3}$ from each other, while any of the integers 0, 3 or 6 (pictured in the lower left triangle of the figure in the middle of Figure \ref{Figure-3-adic}) has 3-adic distance 1 from any of the integers 1, $1+3=4$, $1+2\cdot 3=7$ (pictured in the lower right triangle in figure in the middle of Figure \ref{Figure-3-adic}). 
Hence it is impossible to embed the integers with their 3-adic distance in the Euclidean plane, but Figure \ref{Figure-3-adic} tries to give a visualization of the 3-adic distances of the integers of the form $a_0 + a_1 \cdot 3 + a_2 \cdot 3^2$ with $a_0, a_1, a_2 \in \{0,1\}$ following \cite{Cuoco}. We hope that these figures allow the reader to imagine the fractal-like structure of the 3-adic integers, by adding more and more 3-adic ``digits'' in the presentation  $a_0 + a_1 \cdot 3 + a_2 \cdot 3^2+ a_3 3^3 + \hdots$. In particular, the $p$-adic integers are a compact group whose underlying topological space is totally disconnected. 

As mentioned above, the real numbers $\bR$ and the $p$-adic numbers $\bQ_p$ arise both as completion of the rational numbers $\bQ$, but for a different absolute value. Since the real numbers are connected, more precisely, they are a line, while the $p$-adic numbers are totally disconnected, we have seen that the resulting topological spaces are very different. Another key difference between these fields consist of the following observation: The real numbers contain only one compact subgroup under addition, which is the trivial subgroup consisting only of the element $0$. The $p$-adic numbers, on the other hand, contain infinitely many compact (and open) subgroups under addition, e.g., the $p$-adic integers $\bZ_p$ form a compact, open subgroup of $\bQ_p$, and similarly $p \cdot \bZ_p$ is a compact open subgroup, and more generally $p^n \bZ_p$ for any integer $n$ is a compact, open subgroup of $\bQ_p$. This means in particular that the element 0 in $\bQ_p$ has a basis of open neighborhoods consisting of compact, open subgroups. This property translates into similar properties for $p$-adic groups that shapes the flavor of study of their representation theory.

Before introducing $p$-adic groups, we like to mention that the $p$-adic numbers have a sister field that behaves very similarly to the $p$-adic numbers: the Laurent series $\bF_p((t))$ over a finite field $\bF_p$ with $p$ elements, elements of which can be written as
\[ a_{-n} \cdot t^{-n} + \hdots + a_0  +a_1\cdot t +a_2 \cdot t^2 + \hdots \text{ with }  a_i \in \bF_p , \]
where $t$ is an abstract variable.
There is an obvious similarity between $\bQ_p$ and $\bF_p((t))$ based on the way we presented their elements and the reader is welcome to choose whichever field they feel more comfortable with when reading this article. The key difference between the two fields is that the field of $p$-adic numbers $\bQ_p$ has characteristic zero, it contain the rational numbers $\bQ$ as a subfield, while the field of Laurent series $\bF_p((t))$ has characteristic $p$, it contain the finite field $\bF_p$ with $p$ elements as a subfield in which $p=0$. From now on we write $F$ for any of the two fields $\bQ_p$ or $\bF_p((t))$. While we focus on these two fields in this article, everything covered here also works for finite extensions of these two fields, which are called \textit{non-archimedean local fields}.

In order to introduce $p$-adic groups, we recall that examples of real Lie groups include: 
\begin{itemize}
	\item the general linear group $\GL_n(\bR)$, i.e., $n \times n$-matrices whose determinant is non-zero with the group action being multiplication,
	\item the special linear group $\SL_n(\bR)$, which is the subgroup of $\GL_n(\bR)$ consisting of those matrices whose determinant is 1, 
	\item the special orthogonal group $\SO_n(\bR)$, which is the subgroup of $\SL_n(\bR)$ consisting of those matrices $A$ satisfying $A A^T= \Id$, where $\Id$ denotes the identity matrix and $A^T$ the transpose of $A$, in other words, these are the matrices that preserve the standard inner product, 
	\item the symplectic group $\Sp_{2n}(\bR)$, which is the subgroup of $\GL_{2n}(\bR)$ that preserves a symplectic form.
\end{itemize} 
Replacing the real numbers $\bR$ in the above examples by the field $F$, we obtain examples of what we call \textit{$p$-adic groups}: $\GL_n(F)$, i.e., $\GL_n(\bQ_p)$ and $\GL_n(\bF_p((t)))$, $\SL_n(F)$, $\SO_n(F)$, $\Sp_{2n}(F)$. This means $p$-adic groups are nice subgroups of matrix groups with entries in either the field of $p$-adic numbers $\bQ_p$ or its sister field $\bF_p((t))$. Similar to real Lie groups, there are also $p$-adic groups of type $G_2, F_4, E_6, E_7, E_8$, and these groups are in general classified using combinatorial data (root data) and additional structure related to the absolute Galois group of $F$. Note that we use the terminology ``$p$-adic groups'' for these nice subgroups of $\GL_n(F)$ even though $F$ might be either the $p$-adic numbers $\bQ_p$ or the Laurent series $\bF_p((t))$ over the finite field $\bF_p$. For the experts, the results we discuss below hold for groups $\underline{G}(F)$ where $\underline{G}$ is a connected reductive group over $F$ that splits over a tamely ramified field extension of $F$, but the reader new to this topic is encouraged to just think of the explicit examples mentioned above. For a more detailed introduction to $p$-adic groups including precise definitions we refer the interested reader to Section 2 of the survey article \cite{Fi-CDM}. 

From now on we denote by $G$ a $p$-adic group.

\section{Representations of $p$-adic groups}

\begin{definition}
		A \textit{smooth representation} of $G$ is a group homomorphism 
		\[\pi: G \rightarrow \Aut_\bC(V)\]
		 for some complex vector space $V$ 
		such that for every $v \in V$ there exists an open subgroup $K_v \subset G$ with $\pi(k)(v)=v \, \, \forall \, k \in K_v$. 
\end{definition}
The vector space $V$ underlying a smooth representation of $G$ is usually infinite dimensional. If $G=\SL_2(\bQ_p)$, then we can take as the $K_v$ appearing in the above definition subgroups of the form  $\begin{pmatrix}
	1+p^m\bZ_p & p^m\bZ_p \\
	p^m \bZ_p & 1+p^m\bZ_p
\end{pmatrix}_{\mathrm{det}=1}$ for some positive integers $m$.

Smooth representations of $p$-adic groups have been studied for more than 50 years with a long term goal being to understand the whole category of all smooth representations of $p$-adic groups as explicitly as possible.
\begin{goal}
	Understand the category of all smooth representations of $G$.
\end{goal}

Before explaining the structure of the category of representations of $p$-adic groups, we would like to provide the reader with some examples and a procedure to construct many representations that is used to approach the above long term goal.

\subsection{Example of a representation}\label{section-example-induction}
Let $G=\SL_2(\bQ_p)$ with subgroup $B=\begin{pmatrix}
	\bQ_p & \bQ_p \\ 0 & \bQ_p
\end{pmatrix}_{\mathrm{det}=1}$.

We set 
\begin{align*}
	V := &\left\{ f: B \backslash G \rightarrow \bC \,  \, \mid \, \, f \text{ locally constant} \right\} \\
	\simeq &\left\{ f: 
	\mathbb{P}^1(\bQ_p) \rightarrow \bC \,  \, \mid \, \, f \text{ locally constant}\right\} ,
\end{align*}
which is an infinite dimensional complex vector space, where $\mathbb{P}^1(\bQ_p)$ denotes the projective line over $\bQ_p$.
The map $\pi: G \rightarrow \Aut_\bC(V)$ defining our smooth representation is given as follows: For $g \in G$,  we have
\begin{align*}
	 \pi(g) :& \quad \quad \quad  V \quad \quad \longrightarrow  \quad \quad V \\
	 \pi(g) :&  \quad x \mapsto f(x) \quad  \mapsto  \quad  x \mapsto f(x g)  .
\end{align*}
The requirement that the functions in $V$ are locally constant implies that this is a smooth representation.

\subsection{Parabolic induction}
We would like to generalize the example of the previous subsection. To do so, we introduce the following class of subgroups of $G$ that we will replace the subgroup $B$ above by.
\begin{definition}\label{Def-parabolic}
	A \textit{parabolic subgroup} of $G=\GL_n(F)$ (or $\SL_n(F)$, $\SO_n(F)$ or $\Sp_{2n}(F)$) is a subgroup of the form 
	\[  g \cdot \left(
	\begin{array}{wc{0.01cm}wc{0.01cm}wc{0.01cm}wc{0.01cm}wc{0.01cm}wc{0.01cm}} 
		\cline{1-2}
		\multicolumn{1}{|wc{0.01cm}}{\star}    &    \multicolumn{1}{wc{0.01cm}|}{\star}    & \star    & \star    & \star & \star  \\ 
		\multicolumn{1}{|wc{0.01cm}}{\star}  &     \multicolumn{1}{wc{0.01cm}|}{\star}  & \star    & \star    & \star & \star  \\ \cline{1-3}
		& \multicolumn{1}{wc{0.01cm}|}{}    & \multicolumn{1}{wc{0.01cm}|}{\star}  & \star    & \star & \star  \\ \cline{3-6}
		&  & \multicolumn{1}{wc{0.01cm}|}{} & \star    & \star & \multicolumn{1}{wc{0.01cm}|}{\star}  \\ 
		& \makebox(0,0){\text{\large0}}           &  \multicolumn{1}{wc{0.01cm}|}{}     & \star  & \star & \multicolumn{1}{wc{0.01cm}|}{\star}  \\ 
		&          & \multicolumn{1}{wc{0.01cm}|}{}      &  \star & \star & \multicolumn{1}{wc{0.01cm}|}{\star}  \\ \cline{4-6}
	\end{array}\right) \cdot g^{-1}\] 
	for some choice of number and sizes of blocks (e.g., in the centered expression in this definition we chose three blocks of sizes $2 \times 2, 1\times 1, 3 \times 3$ for $n=6$ to illustrate the shape of the upper-block-triangular matrices that we consider) and some choice of $g \in G$.
\end{definition}
For a more abstract general definition we refer the interested reader to \cite[\S2]{Fi-CDM}.

Any parabolic subgroup admits a decomposition as a semi-direct product of another $p$-adic group, called \textit{Levi subgroup} and denoted $M$ below, corresponding to (the $g$-conjugate of) block-diagonal matrices, and a normal subgroup, called the \textit{unipotent radical} of the parabolic subgroup and denoted by $U$ below, corresponding to (the $g$-conjugate of) the strictly-upper-block-triangular matrices with ones on the diagonal:

\begin{align*}
	\underbrace{ g \left(
	\begin{array}{wc{0.01cm}wc{0.01cm}wc{0.01cm}wc{0.01cm}wc{0.01cm}wc{0.01cm}}
		\cline{1-2}
		\multicolumn{1}{|wc{0.01cm}}{\star}    &    \multicolumn{1}{wc{0.01cm}|}{\star}    & \star    & \star    & \star & \star  \\ 
		\multicolumn{1}{|wc{0.01cm}}{\star}  &     \multicolumn{1}{wc{0.01cm}|}{\star}  & \star    & \star    & \star & \star  \\ \cline{1-3}
		& \multicolumn{1}{wc{0.01cm}|}{}    & \multicolumn{1}{wc{0.01cm}|}{\star}  & \star    & \star & \star  \\ \cline{3-6}
		&  & \multicolumn{1}{wc{0.01cm}|}{} & \star    & \star & \multicolumn{1}{wc{0.01cm}|}{\star}  \\ 
		& \makebox(0,0){\text{\large0}}           &  \multicolumn{1}{wc{0.01cm}|}{}     & \star  & \star & \multicolumn{1}{wc{0.01cm}|}{\star}  \\ 
		&          & \multicolumn{1}{wc{0.01cm}|}{}      &  \star & \star & \multicolumn{1}{wc{0.01cm}|}{\star}  \\ \cline{4-6}
	\end{array}\right) g^{-1}}
=&
\underbrace{
	g
	\left(
	\begin{array}{wc{0.01cm}wc{0.01cm}wc{0.01cm}wc{0.01cm}wc{0.01cm}wc{0.01cm}}
		\cline{1-2}
		\multicolumn{1}{|wc{0.01cm}}{\star}    &    \multicolumn{1}{wc{0.01cm}|}{\star}    &     &     &  &   \\ 
		\multicolumn{1}{|wc{0.01cm}}{\star}  &     \multicolumn{1}{wc{0.01cm}|}{\star}  &     &     & \makebox(0,0){\text{\large0}}  &   \\ \cline{1-3}
		& \multicolumn{1}{wc{0.01cm}|}{}    & \multicolumn{1}{wc{0.01cm}|}{\star}  &     &  &   \\ \cline{3-6}
		&  & \multicolumn{1}{wc{0.01cm}|}{} & \star    & \star & \multicolumn{1}{wc{0.01cm}|}{\star}  \\ 
		& \makebox(0,0){\text{\large0}}           &  \multicolumn{1}{wc{0.01cm}|}{}     & \star  & \star & \multicolumn{1}{wc{0.01cm}|}{\star}  \\ 
		&          & \multicolumn{1}{wc{0.01cm}|}{}      &  \star & \star & \multicolumn{1}{wc{0.01cm}|}{\star}  \\ \cline{4-6}
	\end{array}\right) g^{-1} }
 \cdot
\underbrace{ g  \left(
	\begin{array}{wc{0.01cm}wc{0.01cm}wc{0.01cm}wc{0.01cm}wc{0.01cm}wc{0.01cm}}
		\cline{1-2}
		\multicolumn{1}{|wc{0.01cm}}{1}    &    \multicolumn{1}{wc{0.01cm}|}{0}    & \star    & \star    & \star & \star  \\ 
		\multicolumn{1}{|wc{0.01cm}}{0}  &     \multicolumn{1}{wc{0.01cm}|}{1}  & \star    & \star    & \star & \star  \\ \cline{1-3}
		& \multicolumn{1}{wc{0.01cm}|}{}    & \multicolumn{1}{wc{0.01cm}|}{1}  & \star    & \star & \star  \\ \cline{3-6}
		&  & \multicolumn{1}{wc{0.01cm}|}{} & 1    & 0 & \multicolumn{1}{wc{0.01cm}|}{0}  \\ 
		& \makebox(0,0){\text{\large0}}           &  \multicolumn{1}{wc{0.01cm}|}{}     & 0  &  1& \multicolumn{1}{wc{0.01cm}|}{0}  \\ 
		&          & \multicolumn{1}{wc{0.01cm}|}{}      &  0 & 0 & \multicolumn{1}{wc{0.01cm}|}{1}  \\ \cline{4-6}
	\end{array}\right)  g^{-1}} \\
P   \, \, \qquad \qquad= & \qquad \qquad \, \,  M  \qquad \qquad  \ltimes \qquad \qquad  U
\end{align*}

If $P$ is a proper parabolic subgroup, i.e., $P \subsetneq G$, then the Levi subgroup $M$ is a smaller $p$-adic group. We now introduce a construction to build smooth representations of $G$ from smooth representations of $M$:

\begin{definition}[Parabolic induction]
Let $P= M \ltimes U$ be a parabolic subgroup of $G$. Let $\sigma: M \rightarrow \Aut_\bC(V_\sigma)$ be a smooth representation of $M$. The \textit{parabolic induction} $\Ind_P^G V_\sigma$ is the representation of $G$ with vector space
\[ \Ind_P^G V_\sigma := \left\{ f: G \rightarrow V_\sigma \,   \mid   
\begin{array}{l} 
	f(mug)=\sigma(m)(f(g)) \, \forall \,
	m \in 	M, u \in U, g \in G \! \! \! \! \\
	f \text{ locally constant}
\end{array}
\right\}\]
and for all $g \in G$, 
\begin{align*}
	\Ind_P^G \sigma(g) :& \quad \quad  \Ind_P^G V_\sigma \quad \quad \longrightarrow  \quad \quad \Ind_P^G V_\sigma \\
	\Ind_P^G \sigma(g) :&  \quad \quad x \mapsto f(x) \quad \quad \mapsto  \quad \quad  x \mapsto f(x g)  .
\end{align*}
\end{definition}
Note that we recover the example discussed in Section \ref{section-example-induction} if we set $G=\SL_2(\bQ_p)$, $P=B$ and $V_\sigma=\bC$ with $\sigma: G \mapsto 1 \in \bC^\times=\Aut_\bC(\bC)$.

\subsection{Supercuspidal representations -- the building blocks}

\begin{definition}
	A smooth representation  $\pi: G \rightarrow \Aut_\bC(V)$ is called \textit{irreducible}, if it has exactly two subrepresentations, i.e., if there exist exactly two subspace $W \subseteq V$ (which are $\{0\}$ and $V$) that are preserved under $\pi(g)$ for all $g \in G$.
\end{definition}

\begin{definition}
	An irreducible smooth representation $\pi: G \rightarrow \Aut_\bC(V)$ is called \textit{supercuspidal} if it is not a subrepresentation of
	$\Ind_P^G V_\sigma$
	for every proper parabolic subgroup $P=M\ltimes U \subsetneq G$ and every irreducible smooth representation $\sigma: M \rightarrow \Aut_\bC(V_\sigma)$ of $M$. 
\end{definition}

\begin{fact}
	Let $\pi: G \rightarrow \Aut_\bC(V)$ be an irreducible smooth representation of $G$. Then there exists a parabolic subgroup $P=M\ltimes U \subseteq G$ and a supercuspidal representation $(\sigma, V_\sigma)$ of $M$ such that $V$ is a subrepresentation of $ \Ind_P^G V_\sigma$, which we write as
		\[ V \, \, \hookrightarrow \, \, \Ind_P^G V_\sigma .\]
\end{fact}

Thus the supercuspidal representations form the building blocks of all smooth representations of $p$-adic groups. Therefore our long-term goal to understand all representations of $p$-adic groups can be divided into two steps: constructing all supercuspidal representations, and understanding how the whole category of smooth representations is built up from the supercuspidal representations.

\begin{problem} \label{problem1}
	Construct all supercuspidal representations.
\end{problem}

Solving this first problem of constructing all supercuspidal representations explicitly has plenty of applications within the representation theory of $p$-adic groups, e.g., when studying how representations of $G$ behave when we restrict them to subgroups, or how they behave under various operations done to representations of $p$-adic groups. It is also necessary to describe supercuspidal representations explicitly if one wants to construct an explicit local Langlands correspondence, which means attaching to each representation of a $p$-adic group a number theoretic parameter (roughly a representation of the absolute Galois group of $F$) that bridges between representation theory and number theory. There are also applications to automorphic forms, which are a vast generalization of modular forms, for an example see \cite{FinShin}, and many more applications beyond the representation theory of $p$-adic groups itself.

The vague answer to Problem \ref{problem1} is that we can do this under minor assumptions. We will make this more precise in Section \ref{section-construction-of-sc}.
For now, we first assume that we can construct supercuspidal representations and we try to understand how the whole category of smooth representations of $G$ looks like. 

\subsection{The category of smooth representations}
By Bernstein (\cite{Bernstein}) the category $\Rep(G)$ of all smooth representations of $G$ decomposes into a product of smaller subcategories as follows:
\begin{equation}\label{equation-Bernstein}
	  \underbrace{\mathrm{Rep}(G)}_{\text{smooth} \atop \text{representations}}=\prod_{(M, \sigma)/\sim} \underbrace{\mathrm{Rep}(G)_{[M, \sigma]}}_{\text{Bernstein} \atop \text{block}}  
	\end{equation}
The decomposition is indexed by pairs consisting of a Levi subgroup $M$ of (a parabolic subgroup of) $G$ and a supercuspidal representation $\sigma$ of $M$, where the pairs $(M, \sigma)$ are considered up to an equivalence relation that ensures that each Bernstein block appears only once in the product. More precisely, two pairs $(M, \sigma)$ and $(M', \sigma')$ are equivalent if there exists $g \in G$ and a one-dimensional representation $\phi: M \rightarrow \bC^\times$ that is trivial on all compact subgroups of $M$ such that
\[ M=g M' g^{-1} \quad \text{ and } \quad \sigma(m)=\sigma'(g^{-1}mg)\cdot\phi(m) \, \text{ for } m \in M . \]
Given such a pair $(M, \sigma)$, the attached \textit{Bernstein block} $\mathrm{Rep}(G)_{[M, \sigma]}$ consists of all those representations all of whose irreducible subquotients are contained in some $\Ind_{P'}^G \sigma'$ for some parabolic $P'$ with Levi $M'$ and with $(M', \sigma')$ equivalent to $(M, \sigma)$. In particular, the Bernstein block $\mathrm{Rep}(G)_{[M, \sigma]}$ contains the parabolic induction $\Ind_P^G \sigma$ for any parabolic subgroup $P$ with Levi subgroup $M$, and it can also be characterized as the block that contains $\Ind_P^G \sigma$. Moreover, the Bernstein blocks are indecomposable, i.e., they are not themselves products of smaller non-trivial subcategories, which is why they merit to be called \textit{blocks}.

For example, if $G=\mathrm{SL}_2(\bQ_p)$, then up to conjugation we have two choices for $M$:  
\[ M = G \quad \text{ or } \quad
M=T := \left\{ 	
\left(
\begin{array}{wc{0.2cm}wc{0.2cm}}
\cline{1-1}
\multicolumn{1}{|wc{0.2cm}}{t}    &   	\multicolumn{1}{|wc{0.2cm}}{0}  \\ 
\cline{1-2}
\multicolumn{1}{wc{0.2cm}|}{0}    & \multicolumn{1}{wc{0.2cm}|}{t^{\text{\tiny $-1$}}}   \\
\cline{2-2}
\end{array}
\right)  \mid t \in \bQ_p^\times \right\}. \]

In the case $M=G$, let $\sigma: G \rightarrow \Aut_\bC(V_\sigma)$ be a supercuspidal representation of $G$. Then the Bernstein block $\mathrm{Rep}(\mathrm{SL}_2(\bQ_p))_{[\mathrm{SL}_2(\bQ_p), \sigma]}$ contains as objects arbitrary finite and infinite direct sums of $\sigma$, i.e., $\sigma$, $\sigma \oplus \sigma$, $\sigma \oplus \sigma \oplus \sigma, \hdots$ are objects of this Bernstein block. Since $\sigma$ is irreducible, by Schur's lemma the only morphisms between $V_\sigma$ and $V_\sigma$ that commute with $\sigma(g)$ for all $g \in G$ are multiplication by a scalar, in other words, $\mathrm{Hom}_G(V_\sigma, V_\sigma)= \bC$. From this one can deduce the structure of the morphisms between all other objects in this Bernstein block. 

The Bernstein blocks with more complicated structures arise when one considers $M=T$. As an example, we write $\triv$ for the trivial 
one-dimensional representation of $T$, i.e., $V_{\triv}=\bC$ and $\triv: T \mapsto 1 \in \bC^\times=\Aut_\bC(V_{\triv})$. The Bernstein block 
 $\mathrm{Rep}(\mathrm{SL}_2(\bQ_p))_{[T, \triv]}$ is called the \textit{principal block}  and
contains by definition $\Ind_B^G \triv$, where  $B=\begin{pmatrix}
\bQ_p & \bQ_p \\ 0 & \bQ_p
\end{pmatrix}_{\mathrm{det}=1}$  as above. The representation $\Ind_B^G \triv$ contains the trivial one-dimensional representation $\triv$ of $G$ as a subrepresentation and the quotient $\Ind_B^G \triv / \triv$ is another irreducible smooth representation of $G$ that is called the Steinberg representation and denote by $\mathrm{St}$. Hence the trivial representation and the Steinberg representation are both also contained in the principal block and so is their sum $\triv \oplus \mathrm{St}$. However, the sum $\triv \oplus \mathrm{St}$ is not isomorphic to the parabolic induction $\Ind_B^G \triv$, even though both contain the same subquotients, so the principal block exhibits a more complicated structure.

Nevertheless, it turns out that we can describe this structure rather well as follows: The principal block  $\mathrm{Rep}(\mathrm{SL}_2(\bQ_p))_{[T, \triv]}$ is equivalent to the category of (unital right) modules over the affine Hecke algebra $\cH_{\aff}$ that is defined as follows:
A basis for the complex vector space underlying the affine Hecke algebra $\cH_{\aff}$ is given by 
\[\left\{ T_w \, | \, w \in W_{\text{aff}}:=\langle s_0, s_1 \, | \, s_i^2=1 \rangle \right\}.\]
Note that the abstractly defined group $W_{\text{aff}}$ is the group of transformations of a plane that is generated by two reflections across two parallel hyperplanes.
The relations that describe the multiplication in the affine Hecke algebra are generated by
\[T_{s_i}T_w=\left\{ \, \begin{array}{ll} T_{s_iw} & \ell(s_iw)>\ell(w) \\ pT_{s_iw}+(p-1)T_w & \ell(s_iw)<\ell(w) \end{array} \right. , \]
where $\ell$ denotes the length function $W_{\aff} \rightarrow \bZ_{\geq 0}$ that assigns to an element $w$ in the affine Weyl group $W_{\aff}$ the length of the shortest expression $s_0s_1s_0 \hdots$ or $s_1s_0s_1 \hdots$ that equals $w$.
Finally, the modules of the affine Hecke algebra $\cH_{\aff}$ are well understood.

By the Bernstein decomposition \eqref{equation-Bernstein}, describing the category of smooth representations reduces to the following problem, in addition to Problem \ref{problem1}:
\begin{problem}\label{problem2}
	Describe the Bernstein block $\mathrm{Rep}(G)_{[M, \sigma]}$ explicitly.
\end{problem}

\subsection{Bernstein blocks: Hecke algebras and reduction to depth-zero}

We sketch an answer to Problem \ref{problem2} based on our recent preprints \cite{AFMO1} and \cite{AFMO2} that are joint with Adler, Mishra and Ohara and that were still work in progress during the ECM. For a more detailed description see the survey \cite{Fi-ICM}.

The answer to Problem \ref{problem2} assuming the prime $p$ is not very small (analogous to the assumptions made for the construction of supercuspidal representations discussed in more detail in Section \ref{section-construction-of-sc} below), which we will assume from now on, consists of two parts. 
The first result is of independent interest.

\begin{vaguethm}[{Adler--Fintzen--Mishra--Ohara (\cite{AFMO1} and \cite{AFMO2})}]\label{ThmAFMOA}
Given a Bernstein block $\mathrm{Rep}(G)_{[M, \sigma]}$, we have an equivalence of Bernstein blocks
\[\mathrm{Rep}(G)_{[M, \sigma]} \simeq \mathrm{Rep}(G^0)_{[M^0, \sigma_0]}\]
for some subgroup $G^0 \subseteq G$ (which becomes a Levi subgroup over a field extension), some Levi subgroup $M^0$ of a parabolic subgroup of $G^0$ and some depth-zero supercuspidal representation $\sigma_0$ of $M^0$.
\end{vaguethm}
 We will say a few more words about \textit{depth-zero} supercuspidal representations below, but the main feature of these representations is that they essentially correspond to representations of finite groups of Lie type. Hence the above result allows to reduce a lot of questions in the representation theory of $p$-adic groups and the explicit and categorical local Langlands program to representations of depth zero, where answers are often either already known or much easier to tackle. Such a result has previously only been available for some special cases of $G$ and in the two ``extreme'' cases of $M=T$, where $T$ denotes a split maximal torus, which means in terms of Definition \ref{Def-parabolic} that all the blocks have size $1 \times 1$, by Roche (\cite{Roche98}), and $M=G$ by Ohara (\cite{Ohara-Hecke}).

The second part to answer Problem \ref{problem2} is that we already understood the structure of the Bernstein block $\mathrm{Rep}(G^0)_{[M^0, \sigma_0]}$ very explicitly in many cases since the early 1990s thanks to Morris (\cite{Morris}) and the general case is now treated in  \cite{AFMO1}: We can describe each depth-zero Bernstein block for all $p$-adic groups explicitly as modules over an explicit algebra $\cH$ with explicit generators and relations. More precisely, we proved that the Hecke algebra $\cH$ is isomorphic to a semidirect product of an affine Hecke algebra $\cH_{\aff}$ (generalizing the affine Hecke algebra that we have seen in the example above, allowing more general affine Weyl groups) with a twisted group algebra $\bC[\Omega, \mu]$, see \cite[Notation~3.10.8]{AFMO1} for more precise definitions.
 \begin{vaguethm}[Morris  (\cite{Morris}) and Adler--Fintzen--Mishra--Ohara (\cite{AFMO1})]\label{ThmAFMOB}
 \[ \mathrm{Rep}(G^0)_{[M^0, \sigma_0]} \simeq \bC[\Omega, \mu] \ltimes \cH_{\mathrm{aff}}-\text{mod} \]
 \end{vaguethm}
 
  Thus, combining Theorem \ref{ThmAFMOA} and Theorem \ref{ThmAFMOB} we obtain the following result:
\[ \mathrm{Rep}(G)_{[M, \sigma]} \simeq \mathrm{Rep}(G^0)_{[M^0, \sigma_0]} \simeq \bC[\Omega, \mu] \ltimes \cH_{\mathrm{aff}}-\text{mod}. \]
In fact, the first equivalence (Theorem \ref{ThmAFMOA}) is proven via an explicit isomorphism of Hecke algebras. Since modules over affine Hecke algebras, and more general modules over the semi-direct products  $\bC[\Omega, \mu] \ltimes \cH_{\mathrm{aff}}$ are well studied and understood, this allows us to also understand the structure of the previously mysterious Bernstein blocks.

\subsection{Construction of supercuspidal representations for general $p$-adic groups}\label{section-construction-of-sc}
We return to Problem \ref{problem1}, the construction of (all) supercuspidal representations of $p$-adic groups.

For about 50 years, mathematicians have tried to construct these supercuspidal representations. In 1979, Carayol (\cite{Carayol}) gave a construction of supercuspidal representations that provides all supercuspidal representations of $\GL_n(F)$ for $n$ a prime number different from $p$. In 1986, Moy (\cite{Moy-exhaustion}) proved that Howe's construction (\cite{Howe}) from the 1970s exhausts all supercuspidal representations of $\GL_n(F)$ if $n$ is an integer coprime to $p$. In the early 1990s, Bushnell and Kutzko extended these constructions to obtain all supercuspidal representations of $\GL_n(F)$ for arbitrary $n$ (\cite{BK}).
Similar methods have later been exploited by Stevens (\cite{Stevens}) to construct all supercuspidal representations of classical groups, i.e., symplectic, orthogonal and unitary groups, under the assumption that $p \neq 2$, 
and by S\'echerre and Stevens (\cite{Secherre-Stevens}) to construct all supercuspidal representations of inner forms of $\GL_n(F)$, i.e., $\GL_m(D)$ for some division algebra $D$. There has also been a lot of prior work of mathematicians solving special cases, for a few more details see, for example, \cite{Fi-IHES}.

The picture has been less complete for arbitrary reductive groups.
The introduction of the Moy--Prasad filtration in the 1990s spurred remarkable progress. The work of Moy and Prasad built upon the innovative Bruhat--Tits theory introduced in the 1970s/1980s:
In \cite{BT1, BT2}, Bruhat and Tits defined a building $\cB(G,F)$ associated to $G$ on which $G$ acts. For each point $y$ in $\cB(G,F)$, they constructed a compact subgroup $G_{y,0}$ of $G$, called a \textit{parahoric subgroup}, which is (up to finite index) the stabilizer $G_y$ in $G$ of the point $y$. In \cite{MP1, MP2}, Moy and Prasad defined a filtration of these parahoric subgroups by smaller, normal subgroups
\begin{eqnarray*}
	G_{y,0} \triangleright G_{y, s_1} \triangleright G_{y, s_2} \triangleright G_{y, s_3} \triangleright \hdots  ,
\end{eqnarray*}
where $0< s_1 < s_2 < \hdots$ are real numbers depending on $y$. These subgroups play a crucial role in the study and construction of supercuspidal representations. In some sense, one can heuristically imagine the subgroups $G_{y,s_i}$ as the stabilizer of all points in the Bruhat--Tits building that are within a distance $s_i$ of $y$. We set $G_{y,r}:=G_{y,s_i}$ for $s_{i-1} < r \leq s_i$. 
The quotient $G_{y,0}/G_{y,s_1}$ is a finite group of Lie type.

\begin{figure}[t]
	\includegraphics[width=0.7\textwidth]{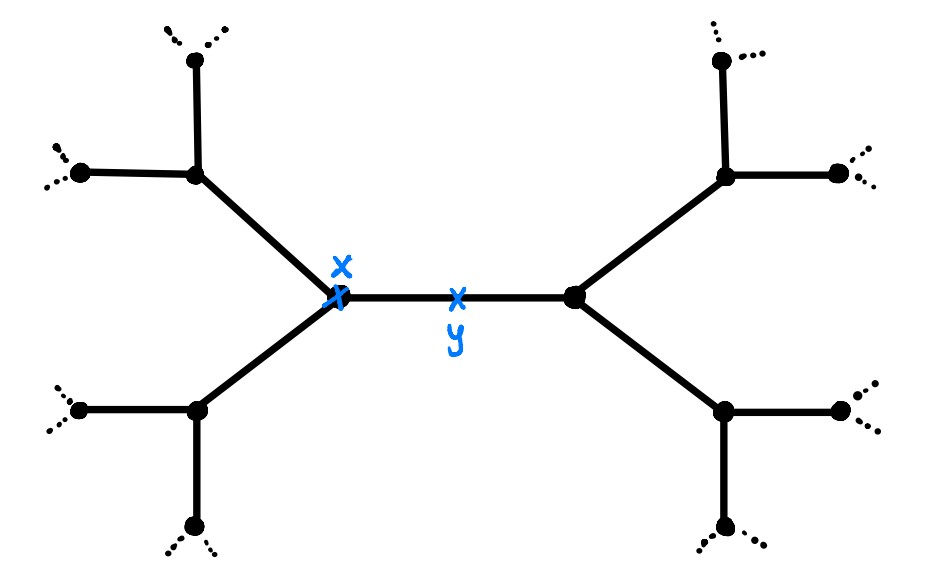}
	\caption{Bruhat--Tits building for $G=\SL_2(\bQ_2)$.}
	\label{Figure-BTtree}
\end{figure}

For example, if we take $G=\SL_2(\bQ_p)$, then the Bruhat--Tits building is an infinite tree with valency $p+1$, see Figure \ref{Figure-BTtree} (for $p=2$). 
Let $y$ be the barycenter of a maximal facet, i.e., the center of an edge of the infinite tree, and let $x$ be a vertex. Then (by choosing an appropriate basis) the associated Moy--Prasad filtrations at the points $x$ and $y$ look like the following: 
\begin{align*}
	G_{x,0} &= \begin{pmatrix} \bZ_p &  \bZ_p \\ \bZ_p &   \bZ_p \end{pmatrix}_{\mathrm{det}=1}  & 	& & G_{y,0}=&\begin{pmatrix} \bZ_p & p \bZ_p \\ \bZ_p &   \bZ_p \end{pmatrix}_{\mathrm{det}=1}   \\
	&   &	 & &   G_{y, 0.5} =&\begin{pmatrix} 1+ p \bZ_p & p \bZ_p \\ \bZ_p &  1+ p \bZ_p \end{pmatrix}_{\mathrm{det}=1}  \\
	G_{x,1} &= \begin{pmatrix} 1+p\bZ_p &  p\bZ_p \\ p\bZ_p &   1+p\bZ_p \end{pmatrix}_{\mathrm{det}=1}  & 	 & &  G_{y, 1} =&\begin{pmatrix} 1+ p \bZ_p & p^2 \bZ_p \\ p \bZ_p &  1+ p \bZ_p \end{pmatrix}_{\mathrm{det}=1} \\
	&  &	&&   G_{y, 1.5} =&\begin{pmatrix} 1+ p^2 \bZ_p & p^2 \bZ_p \\ p \bZ_p &  1+ p^2 \bZ_p \end{pmatrix}_{\mathrm{det}=1}   \\
	G_{x,2} &= \begin{pmatrix} 1+p^2\bZ_p &  p^2\bZ_p \\ p^2\bZ_p &   1+p^2\bZ_p \end{pmatrix}_{\mathrm{det}=1}  & 	 & &  G_{y, 2} =&\begin{pmatrix} 1+ p^2 \bZ_p & p^3 \bZ_p \\ p^2 \bZ_p &  1+ p^2 \bZ_p \end{pmatrix}_{\mathrm{det}=1} \\
	&\vdots   & 	 & &  \vdots & \\
\end{align*}

Based on this filtration, Moy and Prasad introduced in \cite{MP1, MP2} the notion of \textit{depth} of a representation, which measures the first occurrence of a fixed vector in a given representation.  More precisely, it is the smallest non-negative rational number $r$ such that the representation contains a vector fixed under $G_{y,r+}:=\cup_{s>r}G_{y,s}$ for some $y\in \cB(G,F)$. In \cite{MP2}, Moy and Prasad gave a classification of depth-zero representations, showing, roughly speaking, that they correspond to representations of finite groups of Lie type, the group $G_{y,0}/G_{y,0+}$. A similar result was obtained shortly afterwards by Morris (\cite{Morrisdepthzero}) using Hecke algebra techniques.

The first construction of positive-depth supercuspidal representations for general $p$-adic groups\footnote{for experts, as mentioned above, the groups considered are reductive groups over non-archimedean local fields that split over a tamely ramified field extension} was given by Adler (\cite{Adler}) in 1998 and generalized by Yu (\cite{Yu}) in 2001. Since then, Yu's construction has been widely used.
However, about 20 years later it was noticed by Spice that Yu's proof relies on a misprinted (and therefore false) statement in \cite{Gerardin}, and therefore it became uncertain for a few years whether the representations constructed by Yu are irreducible and supercuspidal. In \cite{Fi-Yu-works} we illustrate the impact of this false statement on Yu's proof by providing a counterexample to Proposition 14.1 and Theorem 14.2 of \cite{Yu} and at the same time we offer a different argument for the second half of Yu's proof that avoids the false statements and show that, nevertheless, Yu's construction yields supercuspidal representations.

However, Proposition 14.1 and Theorem 14.2 of \cite{Yu} are the main intertwining results in \cite{Yu} that not only formed the heart of the proof of supercuspidality but are also crucial for other applications. In \cite{FKS}, in joint work with Kaletha and Spice, we therefore construct a quadratic character that allows us to twist Yu's construction so that Yu's initial proof and also the intertwining results hold true for Yu's construction twisted by our quadratic character. This technical result has been the key to calculate Harish-Chandra character formulas (\cite{Spice18} and \cite{Spice21}), to construct a local Langlands correspondence for a large class of representations (\cite{Kaletha-non-singular}), and to obtain Theorems \ref{ThmAFMOA} and \ref{ThmAFMOB} above.

Since many results in our area are only valid for the representations constructed by Yu, one might hope that all supercuspidal representations arise in this way. In 2007 Kim (\cite{Kim}) showed that Yu's representations indeed include all supercuspidal representations if $F$ has characteristic zero, i.e., if we work with $\bQ_p$, but not with $\bF_p((t))$, and if in addition the prime $p$ 
is very large. 
In 2021 we used completely different and much more explicit techniques in \cite{Fi-exhaustion} to show that \textit{all} supercuspidal representations of $G$ are obtained from Yu's construction under much more mild assumptions on $p$, more precisely that $p$ does not divide the order of the absolute Weyl group of $G$, which excludes for exceptional groups at most the primes 2, 3, 5, and 7. In particular, our result holds for any $F$, we do not make a difference between $\bQ_p$ and its sister $\bF_p((t))$.
 Moreover, in \cite{Fi-tame-tori} we illustrate that these assumptions are optimal in general (at least if one considers a slightly more general set-up than we discuss in these notes).
This  settles the long search for all supercuspidal representations under minor assumptions.

On the other hand, this also means that there are more representations to be found for small $p$.
In 2014, Reeder and Yu (\cite{ReederYu}) gave a new construction of some supercuspidal epipelagic representations. Epipelagic representations are representations of smallest positive depth.
In \cite{FR} and \cite{Fi-MP}, we showed that the input for Reeder--Yu's construction also exists for small primes $p$, which provides examples of positive-depth supercuspidal representations that do not arise from Yu's construction. 
These representations have already found various applications to the study of automorphic representations and the inverse Galois problem. 
 However, because these representations are very special and far from exhaustive, many more representations remain to be discovered in the future, a large class of which are the subject of current work in progress. A first big step in this direction is the preprint \cite{FS25} from this year in which we extend in joint work with David Schwein the construction by Yu to also work for $p=2$, a question that had remained open for about 25 years due to the use of the theory of Heisenberg--Weil representations in Yu's construction, which relied on $p \neq 2$. We also remove a superfluous genericity condition by overcoming the challenges of working with disconnected groups, which aligns our construction more closely with potential geometric approaches. This yields new supercuspidal representations of arbitrary depth (for small primes $p$).

All the supercuspidal representations constructed above are of the following form:
\[\text{c-ind}_K^{G} V_\rho := \left\{f: G \rightarrow V_{\rho} \, \left| \, \begin{array}{l} f(kg)=\rho(k)(f(g)) \, \, \forall  g \in G, k \in K \\  f \text{ compactly supported mod center} \end{array} \right. \right\}\]
for some compact-mod-center, open subgroup $K \subset G$, e.g., $\SL_2(\bZ_p) \subset \SL_2(\bQ_p)$, 
and some finite-dimensional representation $(\rho, V_\rho)$ of $K$. The difficulty lies in finding appropriate $K$ and $\rho$ that yield (all) supercuspidal representations.

\begin{ack}
	The author thanks the referee for spotting some typos.
\end{ack}

\begin{funding}
The author was partially supported by the European Research Council (ERC) under the European Union's Horizon 2020 research and innovation programme (grant agreement n° 950326).
\end{funding}

\newpage

\bibliographystyle{alpha}
\bibliography{Fintzenbib}

\end{document}